\documentclass{amsart}

\usepackage{epsfig,color}

\newtheorem{theorem}{Theorem}[section]

\newtheorem{lemma}[theorem]{Lemma}

\theoremstyle{remark}

\newtheorem{remark}[theorem]{Remark}

\newcommand\Q{\mathbb{Q}}

\newcommand\Z{\mathbb{Z}}

\renewcommand\P{\mathbb{P}}

\renewcommand\O{\mathcal{O}}

\newcommand\vx{\overline{x}}
\newcommand\va{\overline{a}}
\newcommand\vb{\overline{b}}
\newcommand\vc{\overline{c}}

\newcounter{nootje}
\begin{document}

\title[Non-Euclidean Pythagorean triples]
      {Non-Euclidean Pythagorean triples, a problem of Euler, and
      rational points on K3 surfaces} 
 \author{Robin Hartshorne}
 \author{Ronald van Luijk}
 \address{Department of Mathematics, Evans Hall, University of
      California, Berkeley, CA 94720-3840}
  \email{robin@math.berkeley.edu}
 \address{Mathematical Sciences Research Institute,
          17 Gauss Way, Berkeley, CA 94720-5070}
  \email{rmluijk@msri.org}

\begin{abstract}
We discover suprising connections between three seemingly different
problems: finding right triangles with rational sides in a non-Euclidean
geometry, finding three integers such that the difference of the
squares of any two is a square, and the problem of finding rational 
points on an algebraic surface in algebraic geometry. We will also 
reinterpret Euler's work on the second problem with a modern point of
view. 
\end{abstract}

\maketitle

\section{Problem I: Pythagorean triples}

An ordinary {\em Pythagorean triple} is a triple $(a,b,c)$ of positive
integers satisfying $a^2+b^2=c^2$. Finding these is equivalent, by
the Pythagorean theorem, to finding right triangles with integral
sides. Since the equation is homogeneous, the problem for rational
numbers is the same, up to a scale factor.

Some Pythagorean triples, such as $(3,4,5)$, have been known since
antiquity. Euclid \cite[X.28, Lemma 1]{euclid}  
gives a method for finding such
triples, which leads to a complete solution of the problem. The
primitive Pythagorean triples are exactly the triples of integers 
$(m^2-n^2,2mn,m^2+n^2)$ for various choices of $m,n$ (up to change of
order). Diophantus in his Arithmetic \cite[Book II, Problem 8]{heath}
mentions the problem of writing any (rational) 
number as the sum of squares. This
inspired Fermat to write his famous ``last theorem'' in the margin.  

Expressed in the language of algebraic geometry, the equation
$x^2+y^2=z^2$ describes a curve in the projective plane. This curve is
parametrized by a projective line according to the assignment 
(in homogeneous coordinates)
\begin{equation}\label{homogparam}
(m,n) \mapsto (m^2-n^2,2mn,m^2+n^2).
\end{equation}
The rational points on the curve correspond to primitive Pythagorean
triples, which explains why the same parametrization appeared above.

Now let us consider the analogous question in a non-Euclidean
geometry. In the hyperbolic plane, if one uses the multiplicative
distance function\footnote{If $A,B$ are two points in the Poincar\'e
  model of a 
  hyperbolic plane, and if $P,Q$ are the ends (on the defining circle)
  of the hyperbolic line containing $A$ and $B$, then $\mu(AB) =
  (AB,PQ)^{-1}$ is the {\em multiplicative distance function} for the
  segment $AB$. Here $(AB,PQ)$ denotes the cross-ratio of the four
  points in the ambient Euclidean plane \cite[39.10]{geb}.}
instead of its logarithm (which is more common),
it makes sense to ask for hyperbolic right triangles whose sides all
have rational numbers $a,b,c$ as lengths. It then follows from the
formulae of hyperbolic trigonometry, see \cite[42.2, 42.3]{geb}, 
that the sines and cosines of the angles of these triangles are 
also rational. The hyperbolic analogue of
the Pythagorean theorem \cite[42.3(f)]{geb} tells us 
$$
\sin \va \cdot \sin \vb = \sin \vc,
$$
where $\va, \vb, \vc$ are the angles of parallelism\footnote{Given a
  segment $AB$ in the hyperbolic plane, let $l$ be the 
hyperbolic line perpendicular to $AB$ at $B$, and let $m$ be the line
through $A$ that is the limiting parallel to $l$ at an end. Then the
angle formed by $AB$ and $m$ is the {\em angle of parallelism} of the
segment $AB$. This construction creates a one-to-one correspondence
(up to congruence) between segments and acute angles in the
hyperbolic plane \cite[40.7.1]{geb}.}
of corresponding segments\footnote{Compare this to the more familiar
  expression of the Pythagorean theorem using the additive distances
  $a,b,c$, namely $\cosh a \cdot \cosh b = \cosh c$.}. 
Since for any segment of length $x$ one has 
\cite[proof of 42.2]{geb} 
$$
\sin \vx = \frac{2x}{1+x^2},
$$
the corresponding arithmetic problem is to find triples $(a,b,c)$ 
of rational numbers, with
\begin{equation}\label{noneuc}
\frac{2a}{1+a^2} \cdot \frac{2b}{1+b^2} = \frac{2c}{1+c^2},
\end{equation}
which we will call {\em non-Euclidean Pythagorean triples}. 
Clearly an even number of the numbers $a,b,c$ are negative, unless 
we have a trivial solution with $abc=0$. After changing signs we may 
assume $a,b,c>0$. Note also that if we replace any number in a non-Euclidean
Pythagorean triple by its inverse, we again get a non-Euclidean
Pythagorean triple, so we may focus on those triples with $a,b,c \geq
1$.  To keep them nontrivial, we require $a,b,c >1$. 
They are not so easy find. (Let the reader try before reading further.)

Note that we must allow rational numbers in stating our problem,
since there are no similar triangles in hyperbolic geometry, or, in
arithmetic terms, since the equation is not homogeneous. 
Indeed, Bjorn
Poonen has shown by an elementary argument that this equation has no
solution in integers $>1$. So our first problem is to find
non-Euclidean Pythagorean triples.

\begin{remark}
%
There are other similar ways to express the hyperbolic Pythagorean
Theorem, for instance \cite[(14.51-52)]{schwerdtfeger} and
\cite{ungar2}. These references do not have an arithmetic point of
view, though. 
\end{remark}

\section{Problem II: A problem of Euler}

Euler, in his Algebra \cite[Part II, \S 236]{alg}, considers the
problem of finding three squares (of integers), $x^2,y^2,z^2$, whose
differences $x^2-y^2$, $x^2-z^2$, $y^2-z^2$ should also be squares.
A first ad hoc argument gives him a single solution 
$(x,y,z)=(697,185,153)$, which turns out to be the smallest
possible, ordered by $|x|$. Then in \S 237 he gives a
method for finding infinitely many solutions. Since one of the
purposes of this paper is to reinterpret Euler's method in terms of
algebraic geometry, we recall his method here. First he notes that 
passing to rational numbers, it is sufficient to find $x,y,z$
satisfying
\begin{equation}\label{diffissq}
\frac{x^2}{z^2} - \frac{y^2}{z^2} = \square, \qquad
\frac{x^2}{z^2} - 1 = \square, \qquad
\frac{y^2}{z^2} - 1 = \square.
\end{equation}
If we set 
\begin{equation}\label{inpq}
\frac{x}{z} = \frac{p^2+1}{p^2-1}, \qquad 
\frac{y}{z} = \frac{q^2+1}{q^2-1}, 
\end{equation}
or $p  = \sqrt{x^2 - z^2}/(x-z)$ and $q = \sqrt{y^2-z^2}/(y-z)$, then 
the second and third equation of (\ref{diffissq}) are automatically
satisfied, just as in the parametrization of Pythagorean triples
mentioned above. Now we only need to satisfy the first
equation. In terms of $p$ and $q$ we want 
$$
\frac{4(p^2q^2-1)(q^2-p^2)}{(p^2-1)^2(q^2-1)^2}
$$
to be a square, and for this it suffices that the numerator be a
square.  Dividing by $4p^2$, we reduce to showing that 
$$
(p^2q^2-1)\left(\frac{q^2}{p^2}-1\right)
$$
is a square. Setting $m=q/p$, we must show that 
\begin{equation}\label{ellp}
(m^2p^4-1)(m^2-1)
\end{equation}
is a square. Obviously, this is a square for $p=1$, but then $x$ is
undefined, so we set $p=1+s$ and then seek to make 
$$
\left(m^2-1 + m^2(s^4+4s^3+6s^2+4s)\right)(m^2-1)
$$
a square. Dividing by $(m^2-1)^2$ and writing $a = m^2/(m^2-1)$ for 
simplicity, we need to make 
\begin{equation}\label{quarts}
1+4as+6as^2+4as^3+as^4
\end{equation}
a square. There are unique $f,g \in \Q[a]$ such
that for $w=1+fs+gs^2$ the coefficients in $w^2$ of $s^k$ for $k=0,1,2$
coincide with the coefficients in (\ref{quarts}), namely $f=2a$ and
$g=3a-2a^2$. Then, to make (\ref{quarts}) equal to $w^2$, we need
\begin{equation}\label{szeroand}
4as^3+as^4 = 2fgs^3+g^2s^4,
\end{equation}
which is the case, besides for $s=0$ (with multiplicity $3$), for 
\begin{equation}\label{sevenbis}
s=\frac{8a-4}{4a^2-8a+1}.
\end{equation}
Reading backwards, take any $m \neq \pm 1$ you like, set $a =
m^2/(m^2-1)$, take $s$ as just given, let $p=1+s$, $q = mp$, and then 
equation (\ref{inpq}) will give $x,y,z$, up to scaling, satisfying the
original problem. Since $m$ is arbitrary, this gives infinitely 
many solutions. 

For example, if we take $m=2$, then 
$$
a= \frac{4}{3}, \quad s = -\frac{60}{23}, \quad p=-\frac{37}{23}, 
\quad q = -\frac{74}{23},
$$
so 
$$
\frac{x}{z} = \frac{949}{420}, \qquad 
\frac{y}{z} = \frac{6005}{4947}, 
$$
giving rise to the relatively prime solution
$$
x = 1564901, \qquad y = 840700, \qquad z = 692580,
$$
and one can easily verify that 
$$
x^2-y^2 = 1319901^2, \qquad
x^2-z^2 = 1403299^2, \qquad
y^2-z^2 = 476560^2.
$$

\begin{remark}\label{equiv}
Euler also considers two other problems. In \S 235 he requires three
integers $a<b<c$, whose sums and differences two at a time are all
squares. One can show easily that this problem is equivalent to ours
as the three pairwise sums of the three integers give a solution to 
our problem and every solution is of this form \cite[section
  4]{leech}. Euler does not mention this equivalence  
and appears not to realize its existence. He does 
see that the squares of $a,b,c$ give a solution to Problem II, but
solutions arising this way tend to be much larger, 
which is why he treats our problem independently by the method
described above. 
%
%
Several other authors have considered these two problems
(see for instance \cite{lyness}, \cite{dickson} and 
references mentioned there). 
Like Euler, however, many did not seem to realize the
equivalence of the two problems. 

In \S 238 Euler requires three squares such that the sum of any two is
again a square. This one has a geometric interpretation in
$3$-dimensional Euclidean space, to find a rectangular box (cuboid)
with integral edges and integral face diagonals. 

Our problem can also be interpreted in terms of cuboids. For if we put 
\begin{equation}\label{cuboids}
x^2-y^2 = t^2, \qquad
x^2-z^2 = u^2, \qquad
y^2-z^2 = v^2,
\end{equation}
these equations are equivalent to 
\begin{equation}\label{fourthquad}
t^2+v^2=u^2, \qquad v^2+z^2=y^2, \qquad t^2+v^2+z^2=x^2.
\end{equation}
Thus our problem is equivalent to finding a cuboid with integral edges
$t,v,z$, of which two face diagonals and the full diagonal are
integral. Note that $(x,u,t)$ is also a solution to Problem II and 
it gives the same cuboid as the triple $(x,y,z)$.

In this connection, it is still an open problem to decide whether or
not there exists a {\em perfect cuboid}, having all edges, face diagonals
and the full diagonal integeral. For more on this problem,
see \cite[D18]{guy}, \cite{leech} and \cite{luijkscr}, and the
references given there. 
\end{remark}

\section{Equivalence of the two problems}

We thank Hendrik Lenstra for first pointing out to us the equivalence
of problems I and II. Equation (\ref{noneuc}) and similar ones
have been studied in relation with Euler's problem before 
\cite[section 4,5]{leech}. Leech also shows how Problem II can be used
to construct spherical right triangles whose sides and angles all have
rational sines and cosines. 

If we take $x,y,z$ satisfying Euler's problem, 
and write equations similar to (\ref{diffissq}), namely
$$
1-\frac{y^2}{x^2} = \square, \qquad
1-\frac{z^2}{x^2} = \square, \qquad
1-\frac{z^2}{y^2}  = \square,
$$
then we can parametrize them (inhomogeneously) as in (\ref{homogparam}) by 
\begin{equation}\label{noneuc2euler}
\frac{y}{x} = \frac{2a}{a^2+1}, \qquad
\frac{z}{x} = \frac{2c}{c^2+1}, \qquad
\frac{z}{y} = \frac{2b}{b^2+1}. 
\end{equation}
From $\frac{y}{x} \cdot \frac{z}{y} = \frac{z}{x}$ we then obtain
equation (\ref{noneuc}) of a non-Euclidean Pythagorean
triple. Conversely, such a 
triple $(a,b,c)$ will give a solution to Euler's problem through
(\ref{noneuc2euler}). 

To find $a,b,c$ explicitly, note that for example the first equation
above, namely 
$$
1-\frac{y^2}{x^2} =\frac{t^2}{x^2},
$$ 
is parametrized by 
$$
\frac{y}{x}= \frac{2a}{a^2+1}, \qquad  \frac{t}{x} = \frac{a^2-1}{a^2+1},
$$
from which we find $a = (x+t)/y$. 
If $x,y,t >0$, then $a > 1$.
%
%
%
%
Similarly we find $b,c>1$. Thus we obtain a one-to-one correspondence between 
solutions $(x,y,z)$ to Problem II with $x>y>z>0$ and $\gcd$ equal to
$1$ on the one hand and (ordered) non-Euclidean Pythagorean triples
$(a,b,c)$ with $a,b,c > 1$ on the other.

So, for example, from Euler's smallest solution $(x,y,z) =
(697,185,153)$ we get
\begin{equation}\label{ex1}
a = \frac{37}{5}, \qquad b= \frac{17}{9}, \qquad c = 9,
\end{equation}
while for the second example above, we obtain
\begin{equation}\label{ex2}
a = \frac{1201}{350}, \qquad b= \frac{97}{51}, \qquad c = \frac{30}{7}.
\end{equation}
It is amusing to verify equation (\ref{noneuc}) for these triples. The
numbers factor, and many of the factors cancel each other as if by magic.


\section{A cycle of five}

In his commentary on the work of Lobachevsky, F. Engel 
noted that to each hyperbolic right triangle, one can associate
another triangle in a natural way. Repeated five times, this process
returns to the original triangle \cite[p. 346-347]{lob}. This
association is
closely related to the formulas of hyperbolic trigonometry, and forms
a parallel to Napier's analogies in spherical trigonometry. 

\begin{figure}[htbp]
\begin{center}
  
\input{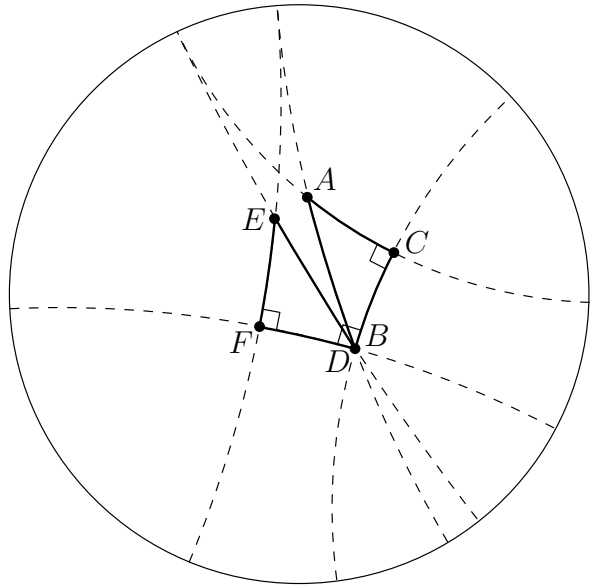}
 
\caption{Engel's associated triangle}
\label{figure:engel}
\end{center}
\end{figure}

Here is the construction. Given the right triangle $ABC$, with sides
$a,b,c$ opposite $A,B,C$ repectively, angles $\alpha, \beta$ at
$A$ and $B$, and a right angle at $C$, draw the perpendicular to $BC$
at $B$, find the limiting parallel to $AB$ that is perpendicular to
this new line, thus obtaining $F$. Draw the limiting parallel from $B$
to $AC$, intersecting the previous limiting parallel at $E$. Then the
new triangle is $DEF$ with $D=B$ (Figure \ref{figure:engel}). 

Note that $\angle FBA$ is complementary to $\beta$, and that it is
the angle of parallelism of segment $DF$. Thus we write $DF=
\overline{\beta'}$, 
where the prime denotes complementary angle, and the bar denotes the
correspondence between segments and angles by the angle of
parallelism. Note also that $\angle EBC$ equals $\overline{a}$, and
so $\angle FDE = \overline{a}'$. Knowing two of the five quantities 
$(d,e,f,\delta,\epsilon)$ of the new triangle, where the variables
denote the obvious lengths and angles, one can compute the
others. Thus, if the original triangle has sides 
and angles $(a,b,c,\alpha,\beta)$, the new triangle has sides and
angles $(b,\overline{\beta'},\overline{\alpha},\overline{a}',
\overline{c})$ (see \cite[42.5 and Exercise 42.23]{geb} for more
details).  It is then easy to verify that this process, repeated five
time, comes back to its starting point. It is an amusing exercise in
hyperbolic trigonometry to compute the new triangle from the old one. 
Since the triangle is determined by any two of its five measurements,
it is enough to compute $b$ (which we already know) and
$\overline{\beta}'$. Here is a recipe. 

\begin{lemma}\label{cycle}
Given $a$ and $b$, set $e= \frac{2a}{a^2-1} \cdot
\frac{b^2-1}{b^2+1}$. Then $\overline{\beta}' = e+\sqrt{1+e^2}$.
\end{lemma}
\begin{proof}
Left to the reader!
\end{proof}

The general formulas for the edges of the new triangle in terms of
those of the old one are not very elegant, but they can be expressed
using a recurrence relation due to Lyness: see \cite{lyness} and 
\cite[p. 524]{leech}. 

Using Lemma \ref{cycle} on an explicit example, say the triple
$$
\left(\frac{37}{5},\frac{17}{9},9\right)
$$
in (\ref{ex1}), we obtain four more triples
\begin{align*}
&\left(\frac{17}{9},\frac{7}{6},\frac{27}{14}\right),\cr
&\left(\frac{7}{6},\frac{5}{4},\frac{21}{16}\right),\cr
&\left(\frac{5}{4},\frac{41}{13},\frac{13}{4}\right), \cr
&\left(\frac{41}{13},\frac{37}{5},13\right)
\end{align*}
as further examples of non-Euclidean Pythagorean triples.
Associated to these are further solutions to Euler's problem, namely
$(x,y,z) = $
\begin{align*}
&(697,185,153), \cr
&(925,765,756), \cr
&(3485,3444,3360), \cr
&(7585,7400,4264), \cr
&(15725, 9061,2405). \cr
\end{align*}
Remarkably, the transformation of order five can be expressed quite
simply in terms of $(x,y,z,t,u,v)$ of (\ref{cuboids}), namely by 
sending $(x,y,z)$ to $(uy,uz,tz)$
and then dividing by the greatest common divisor. 

Just for fun, we computed the $5$-cycle associated to the triple 
in (\ref{ex2}). We get 

\begin{align*}
&\left(\frac{1201}{350},\frac{97}{51},\frac{30}{7}\right),\cr
&\left(\frac{97}{51},\frac{47}{33},\frac{99}{47}\right),\cr
&\left(\frac{47}{33},\frac{37}{23},\frac{1551}{851}\right),\cr
&\left(\frac{37}{23},\frac{73}{26},\frac{74}{23}\right),\cr
&\left(\frac{73}{26},\frac{1201}{350},\frac{40}{7}\right).\cr
\end{align*}

\section{Algebro-geometric interpretation}

The equations (\ref{cuboids})
$$
x^2-y^2 = t^2, \qquad
x^2-z^2 = u^2, \qquad
y^2-z^2 = v^2
$$
describing Euler's problem define a surface $X$ in projective
$5$-space $\P^5$ over the rational numbers. To find integer solutions
to these equations is equivalent to finding rational solutions, hence
to the problem of finding rational points (points with all coordinates
in $\Q$) on the surface $X$. We have thus gone from Euclidean
geometry, through hyperbolic geometry and number theory, back to
geometry, but now algebraic (and arithmetic) geometry. 

For algebraic curves, the problem of finding rational points has been
studied in detail. 
A curve of genus $0$, as soon as it has one rational point, is
isomorphic to $\P^1$, so one knows all its points, and the curve can
be parametrized, just as in the case of the classical Pythagorean
triples. A curve of genus $1$ with at least one rational point is an
elliptic curve. Here one knows that the rational points form a
finitely generated abelian group. It may be finite or infinite. Its
rank can be quite large (Noam Elkies has
found an elliptic curve of rank $28$), but one does not know if there
are elliptic curves of arbitrarily high rank. For curves of genus at
least $2$, Faltings' proof of the Mordell Conjecture tells us that
there are only finitely many rational points. 

For algebraic surfaces (and more generally, varieties of dimension
at least $2$) very little is known about the set of rational points. This
is a topic of intensive active research: see for instance the books 
\cite{pets} and \cite{pots}, and in particular the review paper
\cite{sd} by Swinnerton-Dyer. In our case, the surface $X$
is a complete intersection of three quadric hypersurfaces. 
This implies that outside its singularities the canonical
class $K_X$ is trivial and the irregularity is $0$. Since the only
singularities are ordinary double points, our surface is
(birationally) what is known as a K3 surface. This surface is studied
in detail in \cite{luijkscr}. One knows that there are many K3 surfaces
with no rational points (for example the surface given by
$x_0^4+x_1^4+x_2^4+x_3^4=0$ in $\P^3$, which does not even have any
real points). There are also many whose set of rational points is
dense in the Zariski topology. We will see that ours is one of the
latter. An open problem is whether there exists a K3 surface whose set
of rational points is non-empty, but not dense. Our main result is
the following.

\begin{theorem}\label{main}
The set of rational points on the surface $X$ representing Euler's
problem is dense in the Zariski topology.
\end{theorem}

A good part of the proof is already contained in Euler's calculation. 
The equation 
\begin{equation}\label{affrat}
w^2 = 1+4as+6as^2+4as^3+as^4
\end{equation}
obtained from (\ref{quarts}) describes a surface in the affine
$(a,s,w)$-space. Note that this surface is rational, as its equation
is linear in $a$. By projection it is fibered over the
$a$-line, and for general $a$ the fiber is a curve of genus $1$. Since
it contains some rational points $(s,w) = (0,\pm 1)$, it is an
elliptic curve.  

The assignment $a = m^2/(m^2-1)$ makes the $m$-line a double cover of
the $a$-line. Taking the fibered product of the $m$-line and the
surface given by (\ref{affrat}) over the $a$-line, we get a surface in
$(m,s,w)$-space, fibered by elliptic curves over the $m$-line. This
surface is birational to our surface $X$ via the formulas for
$p,q,x,y,z$ given earlier, and the reverse equations
$$
p = \frac{x+z}{u}, \qquad q = \frac{y+z}{v}.
$$
Thus, our surface is a K3 surface fibered in elliptic curves over a
rational curve, the $m$-line. Its generic fiber is an elliptic curve
$E$ over the rational function field $\Q(m)$, also given by (\ref{affrat}). 
It is worth noting that 
$E$ has $j$-invariant $1728$, making all the elliptic fibers quadratic
twists of the one curve given by $y^2 = x^3-x$. It also implies that
the fibers admit complex multiplication by $\Z[i]$. This is more
clear from (\ref{ellp}), which shows that if $p$ makes (\ref{ellp}) a
square, then so does $ip$. 
Any section of our surface over the $m$-line 
can be given by $\big(m,S(m),W(m)\big)$, where $S(m)$ and
$W(m)$ are rational functions in $m$ that satisfy equation
(\ref{affrat}). The same rational functions determine a point on
$E/\Q(m)$, and thus we obtain a natural
correspondence between sections of our fibration and points on $E$ with
coordinates in $\Q(m)$. Euler's formulas (\ref{sevenbis}) for $s$  and 
$w = 1+fs+gs^2$ give a simple
section of this fibration, thus producing infinitely many rational
points on the surface. 

The new ingredient that we add to this picture is the observation that
after fixing the origin to be $\O=(0,1)$, the $\Q(m)$-point $Q$ that
Euler's section corresponds to, has infinite
order. This can be checked with standard
techniques from the theory of elliptic curves (which we omit).
Taking multiples of this point gives infinitely many more sections
of the fibration. The union of their rational points is dense in the
whole K3 surface, thus proving Theorem \ref{main}.

The point $P = (0,-1)$ has infinite order as well. We claim $Q = 2P$. 
Indeed, note that two points $R$ and $S$ add to $T$ in the group law
on $E$ if and only if there
exists a function $f$ in the function field of the curve whose
principal divisor $(f)$ equals $(T)+(\O)-(R)-(S)$. 
Consider the parabola $F(s,w)=0$ given by $F
= -w+1+2as+(3a-2a^2)s^2$. Euler's calculation shows that this parabola
meets the generic fiber in the point $\O$ with multiplicity $3$ and in
the point $Q$. Since $s$ vanishes at $\O$ and $P$, and the
contribution of the points at infinity to the divisors $(s^2)$ and
$(F)$ are the same, we find that $(F/s^2) = (F)-(s^2)$ equals 
$3(\O)+(Q) - 2(\O)-2(P) = (Q) + (\O) - 2(P)$, so indeed we have $Q =
2P$.

\begin{remark}
By now we have seen several ways of finding rational points on the 
surface $X$. Besides
Euler's first ad hoc argument leading to the smallest integer solution
to Problem II, he found a parametrization of infinitely many 
solutions. His parametrization corresponds to a section of an elliptic
fibration and it turns out to have infinite order, so that there are
in fact infinitely many curves with infinitely many rational points on
the corresponding surface. This shows that the set of rational points 
is dense. Leech \cite[p. 525]{leech} uses a result from Diophantus
(see his Arithmetica, Book V, Lemma 2 to Prop. 7, \cite[p. 205]{heath})
to find a parametric solution corresponding to another curve on the
surface. In terms of algebraic geometry (which is not his viewpoint)
Leech has also found another elliptic fibration. It was already  
mentioned in \cite[page 7]{luijkscr} that Leech's results can be combined
to prove that $X$ contains infinitely many 
parametrizable rational curves, and
thus that the set of rational points on the surface is dense.

There is also a very elementary way to show that there are
infinitely many solutions to Problem II. It comes from
combining the equivalence of Euler's problem and the first problem 
mentioned in Remark \ref{equiv} together with Euler's way of producing a
solution to the former from a solution to the latter, also mentioned
in Remark \ref{equiv}.  Suppose that 
$(x_0,y_0,z_0,t_0,u_0,v_0)$ with $\gcd$ equal to $1$ satisfies
(\ref{cuboids}). Then we have $x_0^2+v_0^2=y_0^2+u_0^2$,
$z_0^2+t_0^2=x_0^2-v_0^2$, 
and $y_0^2-u_0^2 = z_0^2-t_0^2$. This allows us to check that 
$$
(x_1',y_1',z_1',t_1',u_1',v_1') =  (x_0^2+v_0^2,z_0^2+t_0^2,
y_0^2-u_0^2,2x_0v_0,2y_0u_0,2z_0t_0)
$$ 
satisfies (\ref{cuboids}) as well and that $d
=\gcd(x_1',y_1',z_1',t_1',u_1',v_1')$ satisfies $d \leq 2$. 
Note also that we have $x_1'+y_1'+z_1' = x_0^2+y_0^2+z_0^2 \geq 0$.
This means that if we 
start with one solution satisfying $x_0+y_0+z_0>6$, then the new solution 
$(x_1,y_1,z_1,t_1,u_1,v_1)$ we obtain after dividing out $d$ satisfies 
\begin{align*}
x_1+y_1&+z_1 = \frac{1}{d}(x_1'+y_1'+z_1') \geq
\frac{1}{2}(x_1'+y_1'+z_1') \cr
 =& \frac{1}{2}(x_0^2+y_0^2+z_0^2) \geq \frac{1}{6}(x_0+y_0+z_0)^2 >
 x_0+y_0+z_0,
\end{align*}
where the second-to-last inequality is the
quadratic-arithmetic mean inequality. 
We conclude that applying this process repeatedly, obtaining solutions
$(x_n,y_n,z_n)$ for all $n \geq 0$, the expression 
$x_n+y_n+z_n$ increases as $n$ increases, so we get infinitely many 
solutions. We leave the details to the reader. We do not know whether
the set of points arising from a single point by this process can be 
a Zariski dense subset of $X$.
\end{remark}

\begin{remark}
Having found a dense set of rational points on $X$, can we find {\em
  all} rational points on $X$? We believe this is hopeless because a
  special fiber of the fibration may contain many points that 
  do not come from a section of the fibration. 

%
To show this, we first note that the group of rational points 
$E(\Q(m))$ of the generic fiber 
is isomorphic to $\Z \times \Z/2\Z$, generated by 
$P$ and the $2$-torsion point $T=(-2,1)$. This follows from the 
Shioda-Tate formula for the Picard number of our K3 surface
$X$, which equals $20$ by an unpublished result of Beukers and Van
Geemen (see also \cite[p. 8]{luijkscr}). Thus the group of 
$\Q(m)$-rational points on the generic fiber has rank $1$.

Euler finds a parametrization of another rational curve on our surface
that corresponds to a point $R$ on $E$ over
$\Q(\sqrt{a})$, which turns out to satisfy $2R= P+T$. 
This means that for all square values of $a$, we get an extra
point in the fiber above $a$, but the group generated by all points
that Euler finds remains of rank $1$.

Now consider for example Euler's smallest 
solution $(697,185,153)$. It is contained in the fiber that lies above 
$$
m =\frac{q}{p} = \frac{v(x-z)}{u(y-z)} = \frac{13}{5},
$$
or $a = m^2/(m^2-1) = \left(\frac{13}{12}\right)^2$. With a standard technique
from the theory of elliptic curves, namely the so called {\em height
pairing}, implemented in computer algebra packages such as {\sc
  magma}, one can show that in the group law on this fiber the point 
corresponding to Euler's smallest solution is linearly independent of the 
points corresponding to the sections that Euler found, namely $P$,
$Q=2P$, and $R$ satisfying $2R=P+T$. 
The rank over $\Q$ of this special fiber is in
fact $2$. Thus there are many points in special fibers that can not be
obtained as points on sections of 
the surface arising from points in the general fiber. 

Moreover, with a computer we checked that there are $1440$ 
integer solutions $(x,y,z)$ to Euler's Problem II
with $0<z<y<x<10^7$ and $\gcd(x,y,z)=1$. Only $5$ of these  
points lie in a fiber in which that point and the point from Euler's
section $Q$ are linearly dependent. These are the fibers above  
$m \in \left\{2,4,\frac{5}{4},\frac{289}{240},\frac{3267}{2209}\right\}$.

As said before, it therefore seems hopeless that we would ever be 
able to find all rational points. 
In fact, we do not know of any case where the complete set of rational 
points on a K3 surface is described to satisfaction, except 
when there are none. 
\end{remark}

\begin{remark}
Now that we have understood Euler's problem from a modern point of
view, what about the problem of perfect cuboids mentioned above? Here
the coresponding algebraic surface is a surface of general type in
$\P^6$. According to a conjecture of Lang, the rational points on a
surface of general type should be contained in a proper closed
subset. In particular, if there are any (nontrivial points), they will
be much rarer, and difficult to find. So for the moment, this problem
still seems out of reach.
\end{remark}

\end{document}